\documentclass[12pt]{article}

\usepackage{amsmath, epsfig, cite}
\usepackage{amssymb}
\usepackage{amsfonts}
\usepackage{latexsym}
\usepackage{amsthm}

\makeatletter
\renewcommand{\@seccntformat}[1]{{\csname the#1\endcsname}.\hspace{.5em}}
\makeatother

\newtheorem{thm}{Theorem}[section]

\renewcommand{\qed}{\hfill$\Box$\medskip}

\numberwithin{equation}{section}

\begin{document}

\begin{center}
{\large\bf $q$-Analogues of two Ramanujan-type formulas for $1/\pi$}
\end{center}

\vskip 2mm \centerline{Victor J. W. Guo$^1$ and Ji-Cai Liu$^2$\footnote{Corresponding author.}}
\begin{center}
{\footnotesize $^1$School of Mathematical Sciences, Huaiyin Normal
University, Huai'an, Jiangsu 223300,
 People's Republic of China\\
{\tt jwguo@hytc.edu.cn}\\[10pt]

$^2$Department of Mathematics, Wenzhou University, Wenzhou 325035, People's Republic of China\\
{\tt jcliu2016@gmail.com  } }
\end{center}


\vskip 0.7cm \noindent{\bf Abstract.}  We give $q$-analogues of the following two Ramanujan-type formulas for $1/\pi$:
\begin{align*}
\sum_{k=0}^\infty (6k+1)\frac{(\frac{1}{2})_k^3}{k!^3 4^k} =\frac{4}{\pi} \quad\text{and}\quad
\sum_{k=0}^\infty (-1)^k(6k+1)\frac{(\frac{1}{2})_k^3}{k!^3 8^k } =\frac{2\sqrt{2}}{\pi}.
\end{align*}
Our proof is based on two $q$-WZ pairs found by the first author in his earlier work.



\section{Introduction}
In 1997, van Hamme \cite{Hamme} conjectured 13 Ramanujan-type $\pi$ series including
\begin{align}
\sum_{k=0}^\infty (6k+1)\frac{(\frac{1}{2})_k^3}{k!^3 4^k} &=\frac{4}{\pi}, \label{eq:pi-1}\\[5pt]
\sum_{k=0}^\infty (-1)^k(6k+1)\frac{(\frac{1}{2})_k^3}{k!^3 8^k } &=\frac{2\sqrt{2}}{\pi},  \label{eq:pi-2}
\end{align}
have nice $p$-adic analogues (called Ramnujan-type supercongruences). Here we use the Pochhammer symbol $(a)_k=a(a+1)\cdots(a+k-1)$.
All the 13 Ramanujan-type supercongruences have now been confirmed by different authors (see \cite{OZ}). Note that Ekhad and Zeilberger \cite{EZ}
first applied the Wilf--Zeilberger method to prove a Ramanujan-type formula for $\pi$.
Recently, the first author \cite{Guo-1,Guo-2} has formulated $q$-analogues of the (J.2) and (L.2) supercongruences of van Hamme \cite{Hamme},
and confirmed the following special cases:  for any positive odd integer $n$,
\begin{align*}
\sum_{k=0}^{n-1}q^{k^2}[6k+1]\frac{(q;q^2)_k^2 (q^2;q^4)_k }{(q^4;q^4)_k^3}
&\equiv [n](-q)^{\frac{1-n}{2}} \pmod{[n]\Phi_n(q)},
\end{align*}
and, for any odd prime power $n$,
\begin{align*}
\sum_{k=0}^{n-1}(-1)^k [6k+1]\frac{(q;q^2)_k^3}{(q^4;q^4)_k^3}
&\equiv [n](-q)^{-\frac{(n-1)(n+5)}{8}} \pmod{[n]\Phi_n(q)}.
\end{align*}
Here and in what follows, the {\it $q$-shifted factorial} is defined by $(a;q)_n=(1-a)(1-aq)\cdots (1-aq^{n-1})$ for $n\geqslant 1$ and $(a;q)_0=1$,
the {\it $q$-integer} is defined as $[n]=1+q+\cdots+q^{n-1}$, and $\Phi_n(q)$ is the $n$-th {\it cyclotomic polynomial}.

The first purpose of this paper is to prove the following $q$-analogues of \eqref{eq:pi-1} and \eqref{eq:pi-2}, which were originally
conjectured by the first author \cite[Conjecture 4.2]{Guo-1} and \cite[Conjecture 4.5]{Guo-2}, respectively.
\begin{thm}\label{thm:q-pi}
For any complex number $q$ with $|q|<1$, we have
\begin{align}
\sum_{k=0}^{\infty}q^{k^2}[6k+1]\frac{(q;q^2)_k^2(q^2;q^4)_k}{(q^4;q^4)_k^3}
&=\frac{(1+q)(q^2;q^4)_{\infty}(q^6;q^4)_{\infty}}{(q^4;q^4)_{\infty}^2}, \label{a1} \\[5pt]
\sum_{k=0}^{\infty}(-1)^k q^{3k^2}[6k+1]\frac{(q;q^2)_k^3}{(q^4;q^4)_k^3}
&=\frac{(q^3;q^4)_\infty (q^5;q^4)_\infty}{(q^4;q^4)_\infty^2}.  \label{a11}
\end{align}
\end{thm}

To see \eqref{a1} and \eqref{a11} are indeed $q$-analogues of \eqref{eq:pi-1} and \eqref{eq:pi-2},  just notice that the $q$-Gamma function $\Gamma_q(x)$ defined by
$$
\Gamma_q(x)=\frac{(q;q)_\infty}{(q^x;q)_\infty}(1-q)^{1-x},\quad 0<q<1
$$
(see \cite[page 20]{GR}) has the property
$\lim_{q\to 1^{-}}\Gamma_q(x)=\Gamma(x)$, and moreover $\Gamma(x)\Gamma(1-x)=\pi/\sin(\pi x)$.

Z.-W. Sun \cite[(1.6)]{sun-2011} proved that for prime $p\ge 5$,
\begin{align}
\sum_{k=0}^{\frac{p-1}{2}}\frac{{2k\choose k}}{8^k}\equiv \left(\frac{2}{p}\right)
+\left(\frac{-2}{p}\right)\frac{p^2}{4}E_{p-3} \pmod{p^3},  \label{eq:sun}
\end{align}
where $\left(\frac{a}{p}\right)$ is the Legendre symbol modulo $p$ and $E_n$ is the $n$-th Euler number.
The second aim of this paper is to show the following $q$-analogue of \eqref{eq:sun} modulo $p^2$.
\begin{thm}\label{thm:modsun}
For any odd positive integer $n$, we have
\begin{align*}
\sum_{k=0}^{\frac{n-1}{2}}q^{k^2}\frac{(q;q^2)_k}{(q^4;q^4)_k}
\equiv (-q)^{\frac{1-n^2}{8}}\pmod{\Phi_n(q)^2}.
\end{align*}
\end{thm}

We shall prove Theorem \ref{thm:q-pi} in Section 1, and show Theorem \ref{thm:modsun} in
Section 2.

\section{Proof of Theorem \ref{thm:q-pi}}

\noindent{\it Proof of \eqref{a1}.}
We begin with the identity \cite[(2.11)]{Guo-1}:
\begin{align}
\sum_{k=0}^{n-1}q^{k^2}[6k+1]\frac{(q;q^2)_k^2(q^2;q^4)_k}{(q^4;q^4)_k^3}
=\sum_{k=1}^{n}\frac{q^{(n-k)^2}(q^2;q^4)_n(q;q^2)_{n-k}(q;q^2)_{n+k-1}}
{(1-q)(q^4;q^4)_{n-1}^2(q^4;q^4)_{n-k}(q^2;q^4)_{k}}. \label{a2}
\end{align}
For the sake of completeness, we sketch the proof of \cite[(2.11)]{Guo-1} here. Let
\begin{align*}
F(n,k) &=\frac{q^{(n-k)^2}[6n-2k+1](q^2;q^4)_{n}(q;q^2)_{n-k}(q;q^2)_{n+k}}{(q^4;q^4)_{n}^2(q^4;q^4)_{n-k}(q^2;q^4)_k}, \\[5pt]
G(n,k) &=\frac{q^{(n-k)^2}(q^2;q^4)_n (q;q^2)_{n-k}(q;q^2)_{n+k-1}}{(1-q)(q^4;q^4)_{n-1}^2(q^4;q^4)_{n-k}(q^2;q^4)_k},
\end{align*}
where $1/(q^4;q^4)_{m}=0$ for any negative integer $m$. Then
\begin{align}
F(n,k-1)-F(n,k)=G(n+1,k)-G(n,k).  \label{eq:fnk-gnk}
\end{align}
Namely, the functions $F(n,k)$ and $G(n,k)$ form a $q$-WZ pair.  Moreover, the identity \eqref{a2} is equivalent to
\begin{align*}
\sum_{n=0}^{m-1}F(n,0)&= \sum_{k=1}^{m}G\left(m,k\right),
\end{align*}
which follows from \eqref{eq:fnk-gnk} by first summing over $n=0,1,\ldots,m-1$ and then summing over $k$ from $1$ to $m-1$.

Letting $k\to n-k$ on the right-hand side of \eqref{a2}, we obtain
\begin{align}
\sum_{k=0}^{n-1}q^{k^2}[6k+1]\frac{(q;q^2)_k^2(q^2;q^4)_k}{(q^4;q^4)_k^3}
=\sum_{k=0}^{n-1}\frac{q^{k^2}(q^2;q^4)_n(q;q^2)_{k}(q;q^2)_{2n-k-1}}
{(1-q)(q^4;q^4)_{n-1}^2(q^4;q^4)_{k}(q^2;q^4)_{n-k}}. \label{a3}
\end{align}
Furthermore, letting $n\to \infty$ on both sides of \eqref{a3}, we are led to
\begin{align}
\sum_{k=0}^{\infty}q^{k^2}[6k+1]\frac{(q;q^2)_k^2(q^2;q^4)_k}{(q^4;q^4)_k^3}
=\frac{ (q;q^2)_{\infty}}{(q^4;q^4)_{\infty}^2}\sum_{k=0}^{\infty}\frac{q^{k^2}(q;q^2)_{k}}
{(1-q)(q^4;q^4)_{k}}. \label{a4}
\end{align}
Replacing $q$ by $-q$ in Slater's identity \cite[(4)]{Slater-1952}, we have
\begin{align}
\sum_{k=0}^{\infty}q^{k^2}\frac{(q;q^2)_k}{(q^4;q^4)_k}=
\frac{(q^2;q^4)_{\infty}^2}{(q;q^2)_{\infty}}. \label{a5}
\end{align}
The proof of \eqref{a1} then follows from \eqref{a4} and \eqref{a5}.
\qed

\medskip
\noindent{\it Proof of \eqref{a11}.}  We start with the identity \cite[(2.14)]{Guo-2}
\begin{align}
\sum_{k=0}^{n-1}(-1)^k [6k+1]\frac{(q;q^2)_k^3}{(q^4;q^4)_k^3}
&=\sum_{k=1}^{n}\frac{(-1)^{n+k}(q;q^2)_{n+k-1}(q;q^2)_{n-k}^2}{(1-q)(q^4;q^4)_{n-1}^2(q^4;q^4)_{n-k}}, \label{eq:second}
\end{align}
of which the proof is exactly the same as that of \eqref{a2}. The $q$-WZ pair this time is
\begin{align*}
F(n,k) &=\frac{(-1)^{n+k}[6n-2k+1](q;q^2)_{n+k}(q;q^2)_{n-k}^2}{(q^4;q^4)_{n}^2(q^4;q^4)_{n-k}}, \\[5pt]
G(n,k) &=\frac{(-1)^{n+k}(q;q^2)_{n+k-1}(q;q^2)_{n-k}^2}{(1-q)(q^4;q^4)_{n-1}^2(q^4;q^4)_{n-k}}.
\end{align*}

Replacing $q$ by $q^{-1}$ in \eqref{eq:second} and noticing that $(q^{-1};q^{-2})_k^2=(-1)^k q^{k^2}(q;q^2)_k$ and $(q^{-4};q^{-4})_k=(-1)^k q^{4{k+1\choose 2}}(q^{4};q^{4})_k$, we obtain
\begin{align}
\sum_{k=0}^{n-1}(-1)^k q^{3k^2} [6k+1]\frac{(q;q^2)_k^3}{(q^4;q^4)_k^3}
&=\sum_{k=1}^{n}\frac{(-1)^{n+k}q^{(3n+k)(n-k)}(q;q^2)_{n+k-1}(q;q^2)_{n-k}^2}{(1-q)(q^4;q^4)_{n-1}^2(q^4;q^4)_{n-k}} \notag\\[5pt]
&=\sum_{k=0}^{n-1}\frac{(-1)^{k}q^{(4n-k)k}(q;q^2)_{2n-k-1}(q;q^2)_{k}^2}{(1-q)(q^4;q^4)_{n-1}^2(q^4;q^4)_{k}}, \label{eq:second-2}
\end{align}
where the second equality follows from reversing the summation order. Finally, letting $n\to \infty$ on both sides of \eqref{eq:second-2}, we get
\begin{align*}
\sum_{k=0}^{\infty}(-1)^k q^{3k^2}[6k+1]\frac{(q;q^2)_k^3}{(q^4;q^4)_k^3}
&=\frac{(q;q^2)_\infty}{(1-q)(q^4;q^4)_\infty^2},
\end{align*}
since all the summands on the right-hand side of \eqref{eq:second-2} except for the first one ($k=0$) vanish as $n\to\infty$.
The proof of \eqref{a11} then follows from the obvious fact $(q;q^2)_\infty/(1-q)=(q^3;q^4)_\infty(q^5;q^4)_\infty$. \qed

\section{Proof of Theorem \ref{thm:modsun}}
Since
$$
(1-q^{n-2j+1})(1-q^{n+2j-1})+(1-q^{2j-1})^2q^{n-2j+1}=(1-q^n)^2
$$
and $1-q^n\equiv 0\pmod{\Phi_n(q)}$, we have
\begin{align*}
(1-q^{n-2j+1})(1-q^{n+2j-1})\equiv -(1-q^{2j-1})^2q^{n-2j+1}\pmod{\Phi_n(q)^2}.
\end{align*}
Therefore,
\begin{align*}
(-1)^k q^{nk-k^2}(q^{1-n};q^2)_k(q^{n+1};q^2)_k   &=\prod_{j=1}^{k}(1-q^{n-2j+1})(1-q^{n+2j-1})\\[5pt]
&\equiv (-1)^k  \prod_{j=1}^{k}(1-q^{2j-1})^2q^{n-2j+1} \\[5pt]
&=(-1)^kq^{nk-k^2} (q;q^2)_k^2 \pmod{\Phi_n(q)^2}.
\end{align*}
It follows that
\begin{align}
\sum_{k=0}^{\frac{n-1}{2}}q^{k^2}\frac{(q;q^2)_k}{(q^4;q^4)_k}
\equiv \sum_{k=0}^{\frac{n-1}{2}}q^{k^2}\frac{(q^{1-n};q^2)_k(q^{n+1};q^2)_k}{(q;q^2)_k(q^4;q^4)_k}
= (-q)^{\frac{1-n^2}{8}}
 \pmod{\Phi_n(q)^2}. \label{eq:last}
\end{align}
The last identity in \eqref{eq:last} just follows from a terminating $q$-analogue of Whipple's $_3F_2$ sum \cite[Appendix (II.19), page 355]{GR}:
\begin{align*}
_{4}\phi_{3}\left[\begin{array}{c}
q^{-n},\,q^{n+1},\, c,\, -c\\
e,\, c^2q/e,\, -q
\end{array};q,\,q
\right]
=\frac{(eq^{-n},eq^{n+1},c^2q^{1-n}/e,c^2q^{n+2}/e;q^2)_\infty}
{(e,c^2q;q)_\infty} q^{\frac{n(n+1)}{2}}
\end{align*}
with $n\to\frac{n-1}{2}$, $q\to q^2$, $c\to\infty$ and $e\to q$.

\vskip 2mm \noindent{\bf Acknowledgments.} The first author was partially supported by the National Natural Science Foundation of China (grant 11371144), the Natural Science Foundation of Jiangsu Province (grant BK20161304), and the Qing Lan Project of Education Committee of Jiangsu Province.


\begin{thebibliography}{99}
\small \setlength{\itemsep}{-.8mm}

\bibitem{EZ}S.B. Ekhad and D. Zeilberger, A WZ proof of Ramanujan's formula for $\pi$, Geometry,
Analysis, and Mechanics, J.M. Rassias (ed.), World Scientific, Singapore (1994), 107--108.

\bibitem{GR}G. Gasper and M. Rahman, Basic Hypergeometric Series, Second Edition,
Encyclopedia of Mathematics and Its Applications,
Vol. 96, Cambridge University Press, Cambridge, 2004.

\bibitem{Guo-1}V.J.W. Guo, A $q$-Analogue of the (J.2) supercongruence of van Hamme, preprint, 2018.

\bibitem{Guo-2}V.J.W. Guo, A $q$-Analogue of the (L.2) supercongruence of van Hamme, preprint, 2018.

\bibitem{OZ}R. Osburn and W. Zudilin, On the (K.2) supercongruence of Van Hamme,
J. Math. Anal. Appl.  433 (2016),  706--711.

\bibitem{Slater-1952}L.J. Slater, Further identities of the Rogers-Ramanujan type, Proc. London Math. Soc. 54 (1952), 147--167.

\bibitem{sun-2011}Z.-W. Sun, Super congruences and Euler numbers, Sci. China Math. 54 (2011), 2509--2535.

\bibitem{Hamme}L. van Hamme, Some conjectures concerning partial sums of generalized hypergeometric series, $p$-adic functional analysis (Nijmegen, 1996), Lecture Notes in Pure and Appl. Math., vol. 192, Dekker, New York, 1997, 223--236.

\end{thebibliography}
\end{document}